\documentclass{llncs}

\usepackage{amsmath,amssymb,amsfonts, wasysym}
\usepackage[utf8]{inputenc}
\usepackage{hyperref}

\usepackage{amsmath}
\usepackage{url}
\usepackage{listings}

\usepackage{multicol}

\usepackage{graphicx}
\usepackage{caption}
\usepackage{subcaption}

\usepackage{thmtools,thm-restate}

\newcommand{\fixz}{\mathsf{z}}
\newcommand{\pf}{\mathsf{pf}}
\newcommand{\pow}{\mathsf{pow}}
\newcommand{\acc}{\mathsf{acc}}
\newcommand{\eacc}{\mathsf{ech}}

\newcommand{\KNI}[4]{\Phi_{#1, #2,#3}(#4)}
\newcommand{\KNIForm}[1]{\Phi_{#1}}

\newcommand{\PTEL}{\mathbf{PTEL}}
\newcommand{\Psat}{\mathbf{PSAT}}

\newcommand{\Skup}[1]{\mathbf{#1}}
\newcommand{\Fle}{\Skup{For}}
\newcommand{\Sub}[1]{\Skup{Subf}{#1}}
\newcommand{\OverSub}[1]{\Skup{\overline{Subf}}{#1}}

\newcommand{\Nat}{\mathbb{N}}
\newcommand{\SetOfPropLet}{\Skup{Var}}
\newcommand{\SetOfAgents}{\mathbb{A}}
\newcommand{\UnitQ}{[0,1]_{\mathbb{Q}}}
\newcommand{\ClassModels}{\mathtt{Mod}}
\newcommand{\Model}[1]{\mathcal{#1}}
\newcommand{\Powerset}[1]{\mathbb{P}(#1)}
\newcommand{\Ax}{\mathtt{Ax}}

\newcommand{\SetOfAtoms}{\Skup{At}}
\newcommand{\BigLinearsystem}{\Skup{LS}}

\newcommand{\sometime}{\mathtt{F}}
\newcommand{\always}{\mathtt{G}}
\newcommand{\nextT}{\bigcirc}
\newcommand{\prev}{\CIRCLE}
\newcommand{\until}{\mathtt{U}}
\newcommand{\since}{\mathtt{S}}
\newcommand{\once}{\mathtt{P}}
\newcommand{\sofar}{\mathtt{H}}

\newcommand{\know}{\mathtt{K}}
\newcommand{\commonK}{\mathtt{C}}
\newcommand{\groupK}{\mathtt{E}}

\begin{document}
	
	\title{A Probabilistic Temporal Epistemic Logic}
	\author{Zoran Ognjanovi\'{c}\inst{1} \and Angelina Ili\'{c} Stepi\'{c} \inst{1} \and Aleksandar Perovi\'{c}\inst{2}}

\institute{Mathematical Institute of the Serbian Academy of Sciences and Arts, Serbia \\
		\email{[zorano,angelina]@mi.sanu.ac.rs}
		\and
		Faculty of Transport and Traffic Engineering, University of Belgrade, Serbia \\
		\email{pera@sf.bg.ac.rs}
	}

	\maketitle              
	
	\begin{abstract}
We introduce an expressive probabilistic temporal epistemic logic $\PTEL$ suitable to  reason about uncertain knowledge of a non-rigid set of agents that can be changed during time. We define semantics for $\PTEL$ as Kripke models with epistemic accessibility relations for agents' knowledge, linear temporal relations describing possible runs of the system, and probability functions defined on sets of runs. We give an axiomatization of $\PTEL$ and sketch the proof of strong completeness.

	\end{abstract}

\section{Introduction}

Modal-like temporal, epistemic and probabilistic extensions \cite{Hal03,ORM16,Ogn20} of classical logic, as well as their combinations, are used to model different aspects of reasoning in philosophy, artificial intelligence, economics, medicine etc. From the logical point of view interesting questions concern proof--theoretic and modal--theoretic analysis of the corresponding formal systems. Weakly complete axiomatizations (``every valid formula is provable") for several probabilistic logics are given in \cite{FHM90,HM01}, a logic that combine epistemic and probabilistic operators is presented in \cite{FH94}, while \cite{Leh84} analyzes a temporal--epistemic logic and provides a weakly complete axiom system. It turns out that compactness (``a set of formulas is satisfiable iff every finite subset of it is satisfiable") does not hold for these logics, which means that there exists unsatisfiable sets of formulas that are consistent wrt. the corresponding finite axiom systems. This issue can be addressed by providing strongly complete axiomatizations (``every consistent set of formulas is satisfiable"). Such axiom systems (since compactness does not hold) can be obtained by allowing some kind of infiniteness. The papers \cite{OR99,OR00} axiomatize several propositional and first-order probability logics using infinitary rules of inferences, while \cite{Ogn06,MGOS19,TOD20} apply this technique to probabilistic-temporal, temporal-epistemic and probabilistic-epistemic logics. In this approach finiteness is kept at the the object level (the corresponding formal languages are countable, formulas are finite) and infiniteness is present only at the meta level (proofs can be inifinite).

In this paper we introduce an expressive probabilistic temporal multi-agent epistemic logic (denoted $\PTEL$) suitable to reason about uncertain knowledge of a non-rigid set of agents that can be changed during time. We define semantics for $\PTEL$ as Kripke models with epistemic accessibility relations for agents' knowledge, a number of runs consisting of sequences of linearly ordered possible worlds indexed by non-negative integers, and probability functions defined on sets of runs, where runs can be seen as possible executions of a system, so that one can reason about system's properties. Similar models are discussed in \cite{Leh84,FHMV95,Hal03}. The formal language contains epistemic operators for knowledge of a single agent, a group of agents and common knowledge among agents, temporal operators for the past and future in the discrete linear temporal logic, and probabilistic operators that express probabilities of sets of runs. By considering probabilities, $\PTEL$ actually extends the temporal--epistemic logic presented in \cite{MGOS19}.

We give an axiomatization for $\PTEL$ and prove the corresponding strong completeness theorem. In the proofs we will often use induction on the structure of a formula, where the structure is expressed using the rank function $\mathsf{rk}(\cdot)$ which satisfies that:
\begin{enumerate}
\item $\mathsf{rk}(\gamma_1) < \mathsf{rk}(\gamma)$ if $\gamma_1$ is a proper subformula of $\gamma$, and
\item $\mathsf{rk}(\groupK^i \gamma) < \mathsf{rk}(\commonK \gamma)$ for every $i \in \Nat$.
\end{enumerate}
This may be achieved by assigning ordinal ranks to formulas, in particular setting $\mathsf{rk}(\commonK \gamma) := \omega + \mathsf{rk}( \gamma)$. A detailed discussion of such rank function is given in~\cite{BS2009}.


\section{Syntax}
\noindent

Let
\begin{itemize}
  \item $\Nat$ be the set of nonnegative integers,

  \item $\UnitQ$ be the set of all rational numbers from the unit interval,

  \item $\SetOfAgents = \{a_1, \ldots, a_{m}\}$, where $m \in \Nat$, be a set of agents,

  \item $\SetOfPropLet$ be a nonempty at most countable set of propositional letters,

  \item $\Skup{A} = \{A_a | a \in \SetOfAgents\}$ be a subset of $\SetOfPropLet$,

  \item $\Powerset{A}$ be the powerset of the set $A$.
\end{itemize}
The intuitive meaning of the propositional letter $A_a$ is that ``agent $a$ is active''.

The formal language for $\PTEL$ contains the following operators:
\begin{itemize}
  \item classical: $\neg$, $\wedge$,

  \item temporal: $\nextT$, $\until$, $\prev$, $\since$,

  \item epistemic: $\know_a$, $\commonK$, where $a \in \SetOfAgents$,

  \item probabilistic: $\mathtt{P}_{\geqslant s}$, $\mathtt{P}_{a, \geqslant s}$, where $a \in \SetOfAgents$, $s \in \UnitQ$.

\end{itemize}
All operators are unary, except $\wedge$, $\until$ and $\since$ that are binary operators. $\Fle$ denotes the set of formulas defined in the usual way. We will use:
\begin{itemize}
  \item the lowercase Latin letters $p$ and $q$, possibly with indices, to denote propositional variables,

  \item the lowercase Greek letters $\alpha$, $\beta$, $\gamma$, \ldots to denote formulas.
\end{itemize}
A theory is a set of formulas $\Skup{T} \subset \Fle$.

In this paper we will employ the notion of $k$-nested implications.

\begin{definition}[$k$-nested implication]\label{def k nested}
Let $k\in \Nat$.
Let $\Skup{B} = (\beta_{0}, \dots, \beta_{k-1}, \beta_{k})$  be a sequence of $k+1$ formulas, $\alpha \in \Fle$ a formula and $\Skup{X} = (X_{1}, \dots, X_{k-1}, X_{k})$ a sequence of $k$ operators from $ \{ \know_a \, : \, a \in \SetOfAgents\} \cup \{ \nextT, \prev \}$.
The \emph{$k$-nested implication} formula $\KNI{k}{\Skup{B}}{\Skup{X}}{\tau}$ is defined inductively, as follows:
	
	\[
	\KNI{k}{\Skup{B}}{\Skup{X}}{\alpha} =\left\{
	\begin{array}{ll}
	\beta_0 \rightarrow \alpha, \,\, k=0\\
	\beta_{k} \rightarrow X_{k} \KNI{k-1}{\Skup{B}_{j=0}^{k-1}}{\Skup{X}_{j=1}^{k-1}}{\alpha}, \,\, k \geq 1.
	\end{array}
	\right.
	\]
where $\Skup{B}_{j=0}^{k-1} = (\beta_{0}, \dots, \beta_{k-1})$ and $\Skup{X}_{j=1}^{k-1} = (X_{1}, \dots, X_{k-1})$.
\hfill $\blacksquare$
\end{definition}
For example, if $k=4$, $\Skup{X} = (\prev, \know_{a_1}, \nextT, \know_{a_3})$, and $a_1,a_3 \in \SetOfAgents$, then
$$ \KNI{k}{\Skup{B}}{\Skup{X}}{\alpha} = \beta_4 \rightarrow \mathtt{K_{a_3}}( \beta_3 \rightarrow \nextT(\beta_2 \rightarrow \mathtt{K_{a_1}}(\beta_1 \rightarrow \prev(\beta_0 \rightarrow \alpha)))).$$

Also the following abbreviations will be used:
\begin{itemize}

\item if $\Skup{T} \subset \Fle$, then $\nextT \Skup{T} =_{def} \{ \nextT \alpha: \alpha \in \Skup{T} \}$, $\nextT^{i+1} \Skup{T} =_{def} \nextT (\nextT^{i} \Skup{T})$,

\item if $\Skup{T} \subset \Fle$ and $a \in \Skup{A}$, then $\know_a \Skup{T} =_{def} \{ \know_a \alpha: \alpha \in \Skup{T} \}$,

\item if $\Skup{T} \subset \Fle$, then $\Skup{T}^{-\nextT} =_{def} \{ \alpha: \nextT \alpha \in \Skup{T} \}$, $\Skup{T}^{-\nextT^{i+1}} =_{def} (\Skup{T}^{-\nextT^{i}})^{-\nextT}$,

\item if $\Skup{T} \subset \Fle$, then $\Skup{T}^{-\prev} =_{def} \{ \alpha: \prev \alpha \in \Skup{T} \}$, $\Skup{T}^{-\prev^{i+1}} =_{def} (\Skup{T}^{-\prev^{i}})^{-\prev}$,

\item if $\Skup{T} \subset \Fle$ and $a \in \Skup{A}$, then $\Skup{T}^{-\know_a} =_{def} \{ \alpha: \know_a \alpha \in \Skup{T} \}$, $\Skup{T}^{-\know_a^{i+1}} =_{def} (\Skup{T}^{-\know_a^{i}})^{-\know_a}$,

\item $\groupK \alpha =_{def} \bigwedge_{a \in \Skup{A}} \know_a \alpha$, and

\item $\groupK^0 \alpha =_{def} \alpha$ and

\item $\groupK^{n+1} \alpha = \groupK \groupK^{n} \alpha$, $n > 0$

\item $\nextT^0 \alpha =_{def} \alpha; \nextT^{n+1} \alpha=\nextT \nextT^{n} \alpha, n \geqslant 0$,

\item $\prev^0 \alpha =_{def} \alpha; \prev^{n+1} \alpha=\prev \prev^{n} \alpha, n \geqslant 0$.

\item $\sometime \alpha =_{def} (\alpha \rightarrow \alpha) \until \alpha$,

\item $\always \alpha =_{def} \neg \sometime \neg \alpha$,

\item $\once \alpha =_{def} (\alpha \rightarrow \alpha) \since \alpha$,

\item $\sofar \alpha =_{def} \neg \once \neg \alpha$,

\item $\neg \mathtt{P}_{\geqslant s}\alpha =_{def} \mathtt{P}_{<s}\alpha$,
\item $\mathtt{P}_{\geqslant 1-s} \neg \alpha =_{def} \mathtt{P}_{\leqslant s}\alpha$,
\item $\neg \mathtt{P}_{\leqslant s}\alpha =_{def} \mathtt{P}_{> s}\alpha$,
\item $\mathtt{P}_{\geqslant s}\alpha \wedge \mathtt{P}_{\leqslant s}\alpha =_{def} \mathtt{P}_{= s}\alpha$,
\item $\neg \mathtt{P}_{a, \geqslant s}\alpha =_{def} \mathtt{P}_{a, <s}\alpha$,
\item $\mathtt{P}_{a, \geqslant 1-s} \neg \alpha =_{def} \mathtt{P}_{a, \leqslant s}\alpha$,
\item $\neg \mathtt{P}_{a, \leqslant s}\alpha =_{def} \mathtt{P}_{a, > s}\alpha$, and
\item $\mathtt{P}_{a, \geqslant s}\alpha \wedge \mathtt{P}_{a, \leqslant s}\alpha =_{def} \mathtt{P}_{a, = s}\alpha$
\end{itemize}


\section{Semantics}
\noindent

Semantics of $\PTEL$-formulas is given by Kripke-like models with possible worlds that combine temporal, epistemic and probabilistic properties. In this paper we will consider time flow that is isomorphic to the set $\Nat$. The temporal--epistemic part is are essentially the same as models defined in \cite{Leh84,FHMV95}, while \cite{Hal03} discusses how probabilities can be added to this kind of models.

\begin{definition}\label{def_modeli}
A model $\Model{M}$ is any tuple $\langle \Skup{R}, \mathcal{A}, \mathcal{K}, \mathcal{P} \rangle$ such that
\begin{itemize}
		
\item $\Skup{R}$ is a non-empty set of runs, where:
    \begin{itemize}
    \item Every \emph{run} $r$ is a function from $\Nat$ to $\Powerset{\SetOfPropLet}$.
    \item The pair $(r,n)$, where $r \in \Skup{R}$ and $n \in \Nat$, is called a \emph{possible world}; the set of all possible worlds in $\Model{M}$ is denoted by $\Skup{W}$.
    \end{itemize}
		
\item $\mathcal{A}$ is a function from the set of possible world $\Skup{W}$ to $\Powerset{\SetOfAgents}$, where:
    \begin{itemize}
    \item $\mathcal{A}((r,n))$ denotes the set of \emph{active agents} associated to the possible world $(r,n)$, and
    \item $a \in \mathcal{A}((r,n))$ iff $A_a \in r(n)$.
    \end{itemize}
		
\item $\mathcal{K} = \{\mathcal{K}_a: a \in \SetOfAgents\}$ is the set of symmetric and transitive binary \emph{accessibility relations} on $\Skup{W}$, such that:
    \begin{itemize}
    \item $a \not\in \mathcal{A}((r,n))$ iff $(r,n) \mathcal{K}_a (r',n')$ is false for all $(r',n')$.

    \item $\mathcal{K}_a(r,n)$ denotes the set of all possible worlds \emph{accessible}, according to the agent $a$, from $(r,n)$.
    \item If $a \in \mathcal{A}((r,n))$, then
    \begin{itemize}
      \item $(r,n) \mathcal{K}_a (r,n)$.
    \end{itemize}
    \end{itemize}

\item $\mathcal{P}$ is a functions defined on $\Skup{W}$, where $\mathcal{P}((r,n)) = \langle H^{(r,n)}, \mu^{(r,n)}, \{\mathcal{P}_a: a \in \SetOfAgents\}\rangle$ and
    \begin{itemize}
        \item $H^{(r,n)}$ is an algebra of subsets of $\Skup{R}$,
        \item $\mu^{(r,n)}: H^{(r,n)} \rightarrow [0,1]$ is a finitely-additive probability measure on $H^{(r,n)}$,

        \item $\{\mathcal{P}_a: a \in \SetOfAgents\}$ is the set of functions defined on $\Skup{W}$, where $\mathcal{P}_a((r,n)) = \langle \Skup{W}_a^{(r,n)}, H_a^{(r,n)}, \mu_a^{(r,n)} \rangle$ is a probability space such that:
        \begin{itemize}
            \item $\Skup{W}_a^{(r,n)}$ is a non-empty subset of $\Skup{W}$,
            \item $H_a^{(r,n)}$ is an algebra of subsets of $\Skup{W}_a^{(r,n)}$, and
            \item $\mu_a^{(r,n)}: H_a^{(r,n)} \rightarrow [0,1]$ is a finitely-additive probability measure. \hfill $\blacksquare$
        \end{itemize}
    \end{itemize}
\end{itemize}
\end{definition}

Note that:
\begin{itemize}
  \item $r(n)$ is a truth function on propositional letters associated to the possible world $(r,n)$.

  \item The symbols $H^{(r,n)}$ and $\mu^{(r,n)}$ (without the lower indices) correspond to the probabilities defined on runs, while $H_a^{(r,n)}$ and $\mu_a^{(r,n)}$ (with the lower indices) correspond to probabilities related to an agent in a particular possible world and defined on possible worlds.

  \item For every agent $a \in \SetOfAgents$, a particular accessibility relation $\mathcal{K}_a$ in a model, and a particular probability space $\mathcal{P}_a((r,n))$ in every possible world, are defined.

  \item In Definition \ref{def_modeli}, in situations when agents are not active in possible worlds ``dead end worlds'' appear, i.e., worlds from which (according to non active agents) no accessible worlds exist.

  \item If $(r', n')$ is accessible, according to the agent $a$, from $(r,n)$, then $a$ is active in $(r', n')$, i.e., $a \in \mathcal{A}((r',n'))$.

  \item If an agent $a$ is active in all possible worlds, $\mathcal{K}_a$ is an equivalence relation.
\end{itemize}


\section{Satisfiability relation}

\begin{definition} \label{def_satisfiability}
Let $\Model{M} = \langle \Skup{R}, \mathcal{A}, \mathcal{K}, \mathcal{P}\rangle$ be a model. The satisfiability relation $\models$ satisfies:
\begin{enumerate}	
\item if $p \in \SetOfPropLet$, $(r,n) \models p$ iff  $p \in r(n)$,
		
\item $(r,n) \models \alpha \wedge \beta$ iff  $(r,n) \models \alpha$ and $(r,n) \models \beta$,
		
\item $(r,n) \models  \neg \beta$ iff  not $(r,n) \models \beta$ (i.e., $(r,n) \not\models \beta$),
		
\item $(r,n) \models \nextT \beta$ iff $(r, n+1) \models \beta$,

\item $(r,n) \models \alpha \until \beta$ iff there is an integer $j \geqslant n$ such that $(r,j) \models \beta$, and for every integer $k$, such that $n \leqslant k < j$, $(r, k) \models \alpha$,
		
\item $(r,n) \models \prev \beta$ iff $n = 0$, or $n \geqslant 1$ and $(r, n-1) \models \beta$,

\item $(r,n) \models \alpha \since \beta$ iff there is an integer $j \in [0,n] $ such that $(r,j) \models \beta$, and for every integer $k$, such that $j < k \leqslant n$, $(r, k) \models \alpha$,

\item $(r,n) \models \know_a \beta$ iff  $(r',n') \models \beta$ for all $(r',n') \in \mathcal{K}_a(r,n)$,

\item $(r,n) \models \commonK \beta$ iff for every integer $k \geqslant 0$, $(r,n) \models \groupK^{k} \beta$,

\item $(r,n) \models \mathtt{P}_{\geqslant s} \beta$ iff $\mu^{(r,n)}(\{r \in \Skup{R}: (r,0) \models \beta\}) \geqslant s$.
		
\item $(r,n) \models \mathtt{P}_{a,\geqslant s} \beta$ iff $\mu_a^{(r,n)}(\{(r',n') \in \Skup{W}_a^{(r,n)} \, : \, (r',n') \models \beta\}) \geqslant s$. \hfill $\blacksquare$
\end{enumerate}
\end{definition}
Sometimes it is useful to formulate satisfiability of formulas of the form $\commonK \beta$ in an alternative, but equivalent, way. First, the possible world $(r', k')$ is reachable from the possible world $(r, k)$ if there is a finite sequence of possible worlds $(r_0, k_0) = (r, k)$, $(r_1, k_1)$, \ldots, $(r_n, k_n) = (r', k')$ such that for every integer $j \in [0, n-1]$, $(r_j, k_j) \mathcal{K}_{a_j} (r_{j+1}, k_{j+1})$ for some $a_j \in \SetOfAgents$. Then:
\begin{itemize}
\item $(r,n) \models \commonK \beta$ iff for every $(r', k')$ reachable from $(r,n)$, $(r', k') \models \beta$.
\end{itemize}

To ensure that the satisfiability relation is well defined we assume that the sets of the form
\begin{itemize}
  \item $\{r \in \Skup{R} \ : \ (r,0) \models \beta\}$ and
  \item $\{(r',n') \in \Skup{W}_a^{(r,n)} \ : \ (r',n') \models \beta\}$
\end{itemize}
are measurable, i.e., that
\begin{itemize}
  \item $\{r \in \Skup{R} \ : \ (r,0) \models \beta\} \in H^{(r,n)}$ and
  \item $\{(r',n') \in \Skup{W}_a^{(r,n)} \ : \ (r',n') \models \beta\} \in H_a^{(r,n)}$
\end{itemize}
always holds. Models satisfying this condition are measurable. We denote:
\begin{itemize}
  \item the class of all measurable models by $\ClassModels$,

  \item the set $\{r \in \Skup{R} \ : \ (r,0) \models \beta\} \in H^{(r,n)}$ by $[\beta]^{(r,n)}$ and

  \item the set $\{(r',n') \in \Skup{W}_a^{(r,n)} \, : \, (r',n') \models \beta\}$ by $[\beta]_a^{(r,n)}$.
\end{itemize}

\begin{definition} \label{def_validity}
A formula $\alpha$ is valid in a model $\Model{M} \in \ClassModels$ (denoted $\Model{M} \models \alpha$) if for every possible world $(r,n)$ from $\Model{M}$, $(r,n) \models \alpha$.
A formula $\alpha$ is valid (denoted $\models \alpha$) if for every model $\Model{M}$ from $\Model{M} \models \alpha$.

A set of formulas $\Skup{T}$ is satisfied in a possible world $(r,n)$ from a model $\Model{M} \in \ClassModels$ (denoted $(r,n) \models \Skup{T}$), if for every $\alpha \in \Skup{T}$, $(r,n) \models \alpha$. A set of formulas $\Skup{T}$ is satisfiable if there is a possible world $(r,n)$ from a model $\Model{M} \in \ClassModels$ such that $(r,n) \models \Skup{T}$. A formula $\alpha$ is satisfiable if the set $\{ \alpha \}$ is satisfiable.

A formula $\alpha$ is a semantical consequence of the set $\Skup{T}$ of formulas (denoted $\Skup{T} \models \alpha$) if
%
%
%
for every model $\Model{M} \in \ClassModels$ and for every possible world $(r,n)$ from $\Model{M}$, if $(r,n) \models \Skup{T}$, then $(r,n) \models \alpha$.
\hfill $\blacksquare$
\end{definition}
In the above defined notion of semantical consequences,  $\Skup{T} \models \alpha$, formulas from the theory $\Skup{T}$ correspond to local assumptions described in the context of modal logics \cite{Fit93}.

%
%
%


\section{Non-compactness}
\label{noncompactness}

Compactness theorem (``A set of formulas is satisfiable iff all its finite subsets are satisfiable") does not hold for $\PTEL$. Here are some examples of sets of formulas that violate compactness:
\begin{itemize}
  \item $\{\nextT^k \alpha \ : \ k \in \Nat\} \cup \{\neg \always \alpha \}$,
  \item $\{\groupK^k \alpha \ : \ k \in \Nat\} \cup \{\neg \commonK \alpha \}$,
  \item $\{\mathcal{P}_{\leqslant 1/k} \alpha \ : \ k \in \Nat\} \cup \{\neg \mathcal{P}_{= 0} \alpha \}$, etc.
\end{itemize}
Existence of such sets implies that there is no recursive strongly complete axiomatization of our logic \cite{ORM16}. We will provide an infinitary axiomatization such that non-finiteness is present only in meta language:
\begin{itemize}
  \item formulas are finite,
  \item there are inference rules with countably many premisses (and one conclusion), and
  \item proofs are allowed to be infinite.
\end{itemize}


\section{An axiomatization for $\PTEL$}
\label{aksiomatizacija}

The axiomatic system $\Ax$ for $\PTEL$ contains the following axioms and inference rules:

I Propositional axioms and rules\\
%
%
%

\begin{tabbing}
  A$\nextT\prev$C${}_1$. \= $\dfrac{ \{ \KNI{k}{\Skup{B}}{\Skup{X}}{ \neg (( \bigwedge_{l=0}^{i-1} \nextT^l \beta_1) \wedge \nextT^{i} \beta_2)}\, | \,  i \in \Nat\} }{ \KNI{k}{\Skup{B}}{\Skup{X}}{ \neg (\beta_1 \until \beta_2)}}$ \= (Knowledge Necessitation) \kill
   Prop.\>  All instances of classical propositional tautologies \\
   MP. \>  $\dfrac{\alpha, \alpha \to \beta }{\beta }$ \>  (Modus Ponens) \\
\end{tabbing}

II Axioms and rules for reasoning about time\\

\begin{tabbing}
  A$\nextT\prev$C${}_1$. \= $\dfrac{ \{ \KNI{k}{\Skup{B}}{\Skup{X}}{ \neg (( \bigwedge_{l=0}^{i-1} \nextT^l \alpha) \wedge \nextT^{i} \beta)}\, | \,  i \in \Nat\} }{ \KNI{k}{\Skup{B}}{\Skup{X}}{ \neg (\alpha \until \beta)}}$ \= (Knowledge Necessitation) \kill
  A$\nextT\neg$. \>  $\neg \nextT \alpha \leftrightarrow \nextT \neg \alpha$   \\
  A$\nextT\!\!\rightarrow$. \> $\nextT (\alpha \rightarrow \beta) \rightarrow (\nextT \alpha \rightarrow \nextT \beta)$ \> (Distribution Axiom for $\nextT$)  \\
  A$\until\nextT$. \>  $\alpha \until \beta \leftrightarrow \beta \vee (\alpha \wedge \nextT (\alpha \until \beta))$   \\
  A$\until\sometime$. \> $\alpha \until \beta \rightarrow \sometime \beta$   \\
  A$\prev\neg$. \>  $\neg \prev \neg \alpha \rightarrow \prev \alpha $  \\
  A$\prev\!\!\rightarrow$. \>  $\prev (\alpha \rightarrow \beta) \rightarrow (\prev \alpha \rightarrow \prev \beta )$ \> (Distribution Axiom for $\prev$)  \\
  A$\prev\wedge$. \>  $(\prev \alpha \wedge \prev \beta)  \rightarrow \prev (\alpha \wedge \beta)$  \\
  A$\nextT\prev$. \> $\nextT \prev \alpha \leftrightarrow \alpha$ \> (Inversion for $\nextT$ and $\prev$)\\
%
%
%
  A$\nextT\prev$C${}_1$. \>  $\nextT \prev \alpha \rightarrow \prev \nextT \alpha$ \> (Commutativity for $\nextT$ and $\prev$) \\
  A$\nextT\prev$C${}_2$. \> $\neg \prev (\gamma \wedge \neg \gamma) \rightarrow (\nextT \prev \alpha \leftrightarrow \prev \nextT \alpha)$ \> (Commutativity for $\prev$ and $\nextT$) \\
  A$\since\prev$. \>  $\alpha \since \beta \leftrightarrow [\beta \vee (\neg \prev (\alpha \wedge \neg \alpha) \wedge [\alpha \wedge \prev (\alpha \since \beta) ])]$  \\
%
%
  A$\once\prev$. \> $\once \prev \beta$ \\
  R$\nextT$N.  \>  $\dfrac{ \alpha}{ \nextT \alpha }$\label{RNN} \>  (Necessitation for $\nextT$) \\
  R$\prev$N.  \>  $\dfrac{ \alpha}{ \prev \alpha }$\label{RPN} \>  (Necessitation for $\prev$) \\
  R$\until$. \>  $\dfrac{ \{ \KNI{k}{\Skup{B}}{\Skup{X}}{ \neg (( \bigwedge_{l=0}^{i-1} \nextT^l \alpha) \wedge \nextT^{i} \beta)}\, | \,  i \in \Nat\} }{ \KNI{k}{\Skup{B}}{\Skup{X}}{ \neg (\alpha \until \beta)}}$\label{RU} \\
%
  R$\since$. \> $\dfrac{ \{ \KNI{k}{\Skup{B}}{\Skup{X}}{ \neg ( ( \bigwedge_{l=0}^{i-1} \prev^l \alpha) \wedge ( \bigwedge_{l=0}^{i} \neg \prev^l (\alpha \wedge \neg \alpha)) \wedge \prev^{i} \beta)}\, | \,  i \in \Nat\} }{ \KNI{k}{\Skup{B}}{\Skup{X}}{ \neg (\alpha \since \beta)}}$
\end{tabbing}

III Axioms and rules for reasoning about knowledge\\

\begin{tabbing}
  A$\nextT\prev$C${}_1$. \= $\dfrac{ \{ \KNI{k}{\Skup{B}}{\Skup{X}}{ \neg (( \bigwedge_{l=0}^{i-1} \nextT^l \beta_1) \wedge \nextT^{i} \beta_2)}\, | \,  i \in \Nat\} }{ \KNI{k}{\Skup{B}}{\Skup{X}}{ \neg (\beta_1 \until \beta_2)}}$ \= (Knowledge Necessitation) \kill
  A$\know\!\!\rightarrow$. \> $\know_a(\alpha \to \beta) \to  ( \know_a\alpha \to \know_a\beta$) \label{AK} \> (Distribution Axiom for $\know_a$) \\
  A$\know$R. \> $A_a \rightarrow (\know_a \alpha \rightarrow \alpha)$ \> (Reflexivity for $\know_a$) \\
  A$\know$A. \> $A_a \rightarrow \know_a A_a$ \> (Self awareness for $\know_a$) \\
  A$\know$DE. \> $\neg A_a \rightarrow \know_a (\alpha \wedge \neg \alpha)$ \> (Dead end) \\
  A$\know$S. \> $\know_a \neg \alpha \rightarrow \know_a \neg \know_a \alpha$ \> (Symmetry for $\know_a$) \\
  A$\know$T. \> $\know_a \alpha \rightarrow \know_a \know_a \alpha$ \> (Transitivity for $\know_a$) \\
%
%
  A$\commonK\groupK$. \> $\commonK \alpha \to \groupK^m\alpha$, $m \in \Nat$ \label{AC} \\
  R$\know_a$N. \> $\dfrac{ \alpha}{ \know_a \alpha }$\label{RKN} \> (Knowledge Necessitation) \\
  RC. \> $\dfrac{ \{ \KNI{k}{\Skup{B}}{\Skup{X}}{\groupK^i \alpha}\, | \,  i \in \Nat\} }{ \KNI{k}{\Skup{B}}{\Skup{X}}{\commonK \alpha}}$\label{RC}
\end{tabbing}

IV Axioms and rules for reasoning about probability on nruns\\

\begin{tabbing}
  A$\nextT\prev$C${}_1$. \= $\dfrac{ \{ \KNI{k}{\Skup{B}}{\Skup{X}}{ \neg (( \bigwedge_{l=0}^{i-1} \nextT^l \beta_1) \wedge \nextT^{i} \beta_2)}\, | \,  i \in \Nat\} }{ \KNI{k}{\Skup{B}}{\Skup{X}}{ \neg (\beta_1 \until \beta_2)}}$ \= (Knowledge Necessitation) \kill
   AGP1. \>  $\mathtt{P}_{\geqslant 0}{\alpha}$\label{GP1}\>  \\
   AGP2. \>  $\mathtt{P}_{\leqslant r}{\alpha} \to \mathtt{P}_{< t}{\alpha}$, $t>r$\label{GP2}  \\
   AGP3. \>  $\mathtt{P}_{< t}{\alpha} \to \mathtt{P}_{\leqslant t}{\alpha}$\label{GP3}  \\
   AGP4. \>  $(\mathtt{P}_{\geqslant r}{\alpha} \wedge \mathtt{P}_{\geqslant t}{\beta} \wedge \mathtt{P}_{\geqslant 1}{\neg (\alpha \wedge \beta) }) \to \mathtt{P}_{\geqslant \min(1, r+t)}(\alpha \vee \beta) $\label{GP4}  \\
   AGP5. \>  $(\mathtt{P}_{\leqslant r}{\alpha} \wedge \mathtt{P}_{< t}{\alpha}) \to \mathtt{P}_{< r+ t}{(\alpha \vee \beta)} $, $r+t \leqslant 1$\label{GP5}  \\
%
%
   AGP$\prev$. \> $\mathtt{P}_{\geqslant 1}\prev(\alpha \wedge \neg \alpha)$ \label{GP6} \\
   RGPN. \>  $ \dfrac{ \alpha}{ \mathtt{P}_{\geqslant 1}{\alpha} }$  \label{RGP} \>  (Probabilistic Necessitation) \\
   RGA. \>  $ \dfrac{ \{ \KNI{k}{\Skup{B}}{\Skup{X}}{ \mathtt{P}_{\geqslant r - \frac{1}{i}} \alpha} \, | \,  i \geqslant \frac{1}{r} \} }{ \KNI{k}{\Skup{B}}{\Skup{X}}{\mathtt{P}_{\geqslant r} \alpha}}$, $r \in (0,1]_{\mathbb Q}$\label{RGC} \> (Archimedean rule)
\end{tabbing}

V Axioms and rules for reasoning about probability on possible worlds\\

\begin{tabbing}
  A$\nextT\prev$C${}_1$. \= $\dfrac{ \{ \KNI{k}{\Skup{B}}{\Skup{X}}{ \neg (( \bigwedge_{l=0}^{i-1} \nextT^l \beta_1) \wedge \nextT^{i} \beta_2)}\, | \,  i \in \Nat\} }{ \KNI{k}{\Skup{B}}{\Skup{X}}{ \neg (\beta_1 \until \beta_2)}}$ \= (Knowledge Necessitation) \kill
   AP1. \>  $\mathtt{P}_{a,\geqslant 0}{\alpha}$\label{P1}\>  \\
   AP2. \>  $\mathtt{P}_{a,\leqslant r}{\alpha} \to \mathtt{P}_{a,< t}{\alpha}$, $t>r$\label{P2}  \\
   AP3. \>  $\mathtt{P}_{a,< t}{\alpha} \to \mathtt{P}_{a,\leqslant t}{\alpha}$\label{P3}  \\
   AP4. \>  $(\mathtt{P}_{a,\geqslant r}{\alpha} \wedge \mathtt{P}_{a,\geqslant t}{\beta} \wedge \mathtt{P}_{a, \geqslant 1}{\neg (\alpha \wedge \beta) }) \to \mathtt{P}_{a,\geqslant \min(1, r+t)}(\alpha \vee \beta) $\label{P4}  \\
   AP5. \>  $(\mathtt{P}_{a,\leqslant r}{\alpha} \wedge \mathtt{P}_{a,< t}{\alpha}) \to \mathtt{P}_{a,< r+ t}{(\alpha \vee \beta)} $, $r+t \leqslant 1$\label{P5}  \\
   RPN. \>  $ \dfrac{ \alpha}{ \mathtt{P}_{a, \geqslant 1}{\alpha} }$  \label{RP} \>  (Probabilistic Necessitation) \\
   RA. \>  $ \dfrac{ \{ \KNI{k}{\Skup{B}}{\Skup{X}}{ \mathtt{P}_{a,\geqslant r - \frac{1}{i}} \alpha} \, | \,  i \geqslant \frac{1}{r} \} }{ \KNI{k}{\Skup{B}}{\Skup{X}}{\mathtt{P}_{a,\geqslant r} \alpha}}$, $r \in (0,1]_{\mathbb Q}$\label{RA} \> (Archimedean rule)
\end{tabbing}

Note that the axioms from the groups IV and V are similar, which is not surprising since they formalize reasoning about probabilities. But, we need all of them since the language of our logic contains probabilistic operators for both probabilities defined on runs and on possible worlds. The only exception is Axiom AGP$\prev$ in the group IV without a counterpart in the group V. This axiom guarantees that the probability of the set of all runs is $1$, where the set of all runs is defined as the set of all runs beginning with possible worlds in which $\prev(\alpha \wedge \neg \alpha)$ holds.

Using the mentioned abbreviations:
\begin{itemize}
    \item $\mathtt{P}_{\leqslant 1}\alpha = \mathtt{P}_{\geqslant 0} \neg \alpha$
    \item $\mathtt{P}_{a, \leqslant 1}\alpha = \mathtt{P}_{a, \geqslant 0} \neg \alpha$
    \item $\neg \mathtt{P}_{< t}\alpha = \mathtt{P}_{\geqslant t}\alpha$
    \item $\neg \mathtt{P}_{a, < t}\alpha = \mathtt{P}_{a, \geqslant t}\alpha$
    \item $\neg \mathtt{P}_{> t}\alpha = \mathtt{P}_{\leqslant t}\alpha$
    \item $\neg \mathtt{P}_{a, > t}\alpha = \mathtt{P}_{a, \leqslant t}\alpha$
\end{itemize}
we have the following versions of the above axioms:
\begin{tabbing}
  A$\nextT\prev$C${}_1$. \= $\dfrac{ \{ \KNI{k}{\Skup{B}}{\Skup{X}}{ \neg (( \bigwedge_{l=0}^{i-1} \nextT^l \beta_1) \wedge \nextT^{i} \beta_2)}\, | \,  i \in \Nat\} }{ \KNI{k}{\Skup{B}}{\Skup{X}}{ \neg (\beta_1 \until \beta_2)}}$ \= (Knowledge Necessitation) \kill
   AGP1'. \>  $\mathtt{P}_{\leqslant 1}{\alpha}$ \>  \\
   AGP2'. \>  $\mathtt{P}_{\geqslant t}{\alpha} \to \mathtt{P}_{> r}{\alpha}$, $t>r$  \\
   AGP3'. \>  $\mathtt{P}_{> t}{\alpha} \to \mathtt{P}_{\geqslant t}{\alpha}$  \\
   AP1'. \>  $\mathtt{P}_{a, \leqslant 1}{\alpha}$ \>  \\
   AP2'. \>  $\mathtt{P}_{a, \geqslant t}{\alpha} \to \mathtt{P}_{a, > r}{\alpha}$, $t>r$  \\
   AP3'. \>  $\mathtt{P}_{a, > t}{\alpha} \to \mathtt{P}_{a, \geqslant t}{\alpha}$
\end{tabbing}

The rules R$\until$, R$\since$, RC, RGA and RA are infinitary, i.e., each of them has a countable set of assumptions and one conclusion. An equivalent form of Rule RGA is:
\begin{tabbing}
  A$\nextT\prev$C${}_1$. \= $\dfrac{ \{ \KNI{k}{\Skup{B}}{\Skup{X}}{ \neg (( \bigwedge_{l=0}^{i-1} \nextT^l \beta_1) \wedge \nextT^{i} \beta_2)}\, | \,  i \in \Nat\} }{ \KNI{k}{\Skup{B}}{\Skup{X}}{ \neg (\beta_1 \until \beta_2)}}$ \= (Knowledge Necessitation) \kill
   RGA'. \>  $ \dfrac{ \{ \KNI{k}{\Skup{B}}{\Skup{X}}{ \mathtt{P}_{\leqslant r + \frac{1}{i}} \alpha} \, | \,  i \in \Nat \} }{ \KNI{k}{\Skup{B}}{\Skup{X}}{\mathtt{P}_{\leqslant r} \alpha}}$.
\end{tabbing}
A similar form can be given for Rule RA.


\begin{definition} \label{def_prof}
Let $\lambda$ be a finite or countable ordinal.

A formula $\alpha $ is a \emph{theorem}, denoted by $\vdash \alpha$, if there is an at most countable sequence of formulas $\alpha_0$, $\alpha_1$, \ldots , $\alpha_{\lambda+1}$ from $\Fle$, such that:
\begin{itemize}
  \item $\alpha_{\lambda+1}=\alpha$, and
  \item every $\alpha_i$ is an instance of some axiom schemata or is obtained from the preceding formulas by an application of an inference rule.
\end{itemize}
A formula $\alpha$ is \emph{derivable from} a set $\Skup{T}$ of formulas ($\Skup{T} \vdash \alpha$) if there is an at most countable sequence of formulas $\alpha_0$, $\alpha_1$, \ldots , $\alpha_{\lambda+1}$ from $\Fle$ such that:
\begin{itemize}
  \item $\alpha_{\lambda+1}=\alpha$, and
  \item every $\alpha_i$ is an instance of some axiom schemata or a formula from the set $\Skup{T}$, or it is obtained from the previous formulas by an inference rule, with the exception that the premises of the inference rules R$\nextT$N, R$\prev$N, R$\know_a$N, RGPN and RPN must be theorems .
\end{itemize}
The corresponding sequence of formulas is a \textit{proof} for $\alpha$ (from the set $\Skup{T}$). \hfill $\blacksquare$
\end{definition}
Note that the length of a proof is an \emph{at most countable successor ordinal}, and that a formula is a theorem iff it is derivable from the empty set.

\begin{definition}\label{def_consist}
A set $\Skup{T}$ of formulas is \emph{inconsistent} w.r.t. $\Ax$ (or simply \emph{inconsistent}) if $\Skup{T} \vdash \alpha$ for every formula $\alpha$, otherwise it is $\Ax$-\emph{consistent} (or simply \emph{consistent}).  A set $\Skup{T}$ of formulas is $\Ax$-\emph{maximal consistent} (or simply \emph{maximal consistent}) if it is consistent, and each proper superset of $T$ is inconsistent.

A set of formulas $\Skup{T}$ is \emph{deductively closed} w.r.t. $\Ax$ (or simply \emph{deductively closed}) if it contains all the formulas derivable from $\Skup{T}$, i.e.,  $\alpha \in \Skup{T}$ whenever $\Skup{T} \vdash \alpha$. \hfill $\blacksquare$
\end{definition}

The infinitary inference rules in the above given system $\Ax$ guarantee that sets from Section \ref{noncompactness} that violate compactness are inconsistent. The next example illustrates this property for the set with the probabilistic operators.

\begin{example}
Let $\Skup{T} = \{\mathtt{P}_{\leqslant 1/k} \alpha \ : \ k \in \Nat\} \cup \{\neg \mathtt{P}_{= 0} \alpha \}$. Then:
\begin{itemize}
  \item[] $\Skup{T} \vdash \neg \mathtt{P}_{=0} \alpha$, since $\neg \mathtt{P}_{=0} \alpha \in \Skup{T}$
  \item[] $\Skup{T} \vdash \mathtt{P}_{\leqslant 1/k} \alpha$, since $\mathtt{P}_{\leqslant 1/k} \alpha \in \Skup{T}$, for every $k \in \Nat$
  \item[] $\Skup{T} \vdash \mathtt{P}_{\leqslant 0} \alpha$, by Rule RGA'
  \item[] $\Skup{T} \vdash \mathtt{P}_{= 0} \alpha$, by Axiom AGP1 and definition of $\mathtt{P}_{=0}$
  \item[] $\Skup{T} \vdash \bot$.\hfill $\blacksquare$
\end{itemize}
\end{example}


\section{Auxiliary statements} 

\begin{theorem}[Deduction theorem]\label{deduction}
If $\Skup{T} \subset \Fle$, then
$$\Skup{T}, \{ \alpha \} \vdash \beta \mbox{  iff  } \Skup{T} \vdash \alpha \rightarrow \beta.$$
\end{theorem}
\begin{proof}
The $(\leftarrow)$-direction is standard.
To prove the $(\rightarrow)$-direction we use transfinite induction on the length of the proof of $\beta$ from $\Skup{T} \cup \{ \alpha \}$. The next case for Rule R$\since$ illustrates the idea, while the cases for other rules can be proved analogously.
Let us assume $\beta =  \KNI{k}{\Skup{B}}{\Skup{X}}{ \neg (\gamma \since \delta)}$, and:
\begin{itemize}

\item[] $\Skup{T}, \alpha \vdash \KNI{k}{\Skup{B}}{\Skup{X}}{ \neg ( ( \bigwedge_{l=0}^{i-1} \prev^l \gamma) \wedge ( \bigwedge_{l=0}^{i} \neg \prev^l (\gamma \wedge \neg \gamma)) \wedge \prev^{i} \delta)}$, for every $i \in \Nat$.

\end{itemize}
If $k>0$, then:
\begin{itemize}
\item[] $\Skup{T} \vdash \alpha \rightarrow \KNI{k}{\Skup{B}}{\Skup{X}}{ \neg ( ( \bigwedge_{l=0}^{i-1} \prev^l \gamma) \wedge ( \bigwedge_{l=0}^{i} \neg \prev^l (\gamma \wedge \neg \gamma)) \wedge \prev^{i} \delta)}$, for every $i \in \Nat$, by the induction hypothesis

\item[] $\Skup{T} \vdash \alpha \rightarrow (\beta_{k} \rightarrow X_{k} \KNI{k-1}{\Skup{B}_{j=0}^{k-1}}{\Skup{X}_{j=1}^{k-1}}{ \neg ( ( \bigwedge_{l=0}^{i-1} \prev^l \gamma) \wedge ( \bigwedge_{l=0}^{i} \neg \prev^l (\gamma \wedge \neg \gamma)) \wedge \prev^{i} \delta)}$, for every $i \in \Nat$, by the definition of $\KNIForm{k}$

\item[] $\Skup{T} \vdash (\alpha \wedge \beta_{k}) \rightarrow X_{k} \KNI{k-1}{\Skup{B}_{j=0}^{k-1}}{\Skup{X}_{j=1}^{k-1}}{ \neg ( ( \bigwedge_{l=0}^{i-1} \prev^l \gamma) \wedge ( \bigwedge_{l=0}^{i} \neg \prev^l (\gamma \wedge \neg \gamma)) \wedge \prev^{i} \delta)}$, for every $i \in \Nat$, by propositional reasoning.
\end{itemize}
Now, we define $\Skup{\overline{B}} = (\beta_{0}, \dots, \beta_{k-1}, \alpha \wedge \beta_{k})$, and have:
\begin{itemize}
\item[] $\Skup{T} \vdash \KNI{k}{\Skup{\overline{B}}}{\Skup{X}}{ \neg ( ( \bigwedge_{l=0}^{i-1} \prev^l \gamma) \wedge ( \bigwedge_{l=0}^{i} \neg \prev^l (\gamma \wedge \neg \gamma)) \wedge \prev^{i} \delta)}$, for every $i \in \Nat$

\item[] $\Skup{T} \vdash \KNI{k}{\Skup{\overline{B}}}{\Skup{X}}{ \neg (\gamma \since \delta)}$, by R$\since$

\item[] $\Skup{T} \vdash (\alpha \wedge \beta_{k}) \rightarrow X_{k} \KNI{k-1}{\Skup{B}_{j=0}^{k-1}}{\Skup{X}_{j=1}^{k-1}}{ \neg (\gamma \since \delta)}$

\item[] $\Skup{T} \vdash \alpha \rightarrow (\beta_{k} \rightarrow X_{k} \KNI{k-1}{\Skup{B}_{j=0}^{k-1}}{\Skup{X}_{j=1}^{k-1}}{ \neg (\gamma \since \delta)}$

\item[] $\Skup{T} \vdash \alpha \rightarrow \KNI{k}{\Skup{B}}{\Skup{X}}{ \neg (\gamma \since \delta)}$.
\end{itemize}
If $k=0$, we reason as above, with the proviso that $\Skup{\overline{B}} = (\alpha \wedge \beta_{0})$.
\hfill $\blacksquare$
\end{proof}


\begin{theorem}[Strong necessitation]\label{strongnec}
If $\Skup{T} \subset \Fle$ and $\Skup{T} \vdash \gamma$, then
\begin{enumerate}
  \item $\nextT \Skup{T} \vdash \nextT \gamma$, \label{strongnec_next}
  \item $\prev \Skup{T} \vdash \prev \gamma$, and \label{strongnec_prev}
  \item $\know_a \Skup{T} \vdash \know_a \gamma$, for every $a \in \Skup{A}$. \label{strongnec_k}
\end{enumerate}
\end{theorem}
\begin{proof}
We use transfinite induction on the length of the proof of $\gamma$ from $\Skup{T}$. Recall that $\Skup{T} \vdash \gamma$ means that there is an at most countable sequence of formulas $\gamma_0$, $\gamma_1$, \ldots , $\gamma_{\lambda+1}$ such that:
\begin{itemize}
  \item $\lambda$ is a finite or countable ordinal,
  \item $\gamma_{\lambda+1}=\gamma$, and
  \item every $\gamma_i$ is:
   \begin{itemize}
     \item an instance of some axiom schemata,
     \item a formula from the set $\Skup{T}$, or
     \item obtained from the previous formulas by an inference rule (note that premises of R$\nextT$N, R$\prev$N, R$\know_a$N, RPN and must be theorems).
   \end{itemize}
\end{itemize}

\noindent
Here we consider the above statement (\ref{strongnec_next}):
\begin{itemize}
  \item if $\Skup{T} \vdash \gamma$, then $\nextT \Skup{T} \vdash \nextT \gamma$.
\end{itemize}
If $\gamma_i$ is an instance of an axiom schemata, $\vdash \gamma_i$, then also $\vdash \nextT \gamma_i$, and obviously $\nextT \Skup{T} \vdash \nextT \gamma_i$. If $\gamma_i \in \Skup{T}$, then $\nextT \gamma_i \in \nextT \Skup{T}$, and $\nextT \Skup{T} \vdash \nextT \gamma_i$. To illustrate the idea of the proof let us assume that $\gamma_i$ is obtained by one of the infinitary rules, e.g. by the rule R$\until$. Then:
\begin{itemize}
  \item[] $\Skup{T} \vdash \KNI{k}{\Skup{B}}{\Skup{X}}{ \neg (( \bigwedge_{l=0}^{i-1} \nextT^l \alpha) \wedge \nextT^{i} \beta)}$, for every $i \in \Nat$

  \item[] $\Skup{T} \vdash \KNI{k}{\Skup{B}}{\Skup{X}}{ \neg (\alpha \until \beta)}$, by R$\until$.
\end{itemize}
Then, we have:
\begin{itemize}
  \item[] $\Skup{T} \vdash \KNI{k}{\Skup{B}}{\Skup{X}}{ \neg (( \bigwedge_{l=0}^{i-1} \nextT^l \alpha) \wedge \nextT^{i} \beta)}$, for every $i \in \Nat$

  \item[] $\nextT \Skup{T} \vdash \nextT \KNI{k}{\Skup{B}}{\Skup{X}}{ \neg (( \bigwedge_{l=0}^{i-1} \nextT^l \alpha) \wedge \nextT^{i} \beta)}$, for every $i \in \Nat$, by the induction hypothesis, and

  \item[] $\nextT \Skup{T} \vdash (\beta_{k} \vee \neg \beta_{k}) \rightarrow \nextT \KNI{k}{\Skup{B}}{\Skup{X}}{ \neg (( \bigwedge_{l=0}^{i-1} \nextT^l \alpha) \wedge \nextT^{i} \beta)}$, for every $i \in \Nat$.
\end{itemize}
We extend $\Skup{B}$ and $\Skup{X}$:
\begin{itemize}
  \item $\Skup{\overline{B}} = (\beta_{0}, \dots, \beta_{k}, \beta_{k} \vee \neg \beta_{k})$, and
  \item $\Skup{\overline{X}} = (X_{1}, \dots, X_{k}, \nextT)$
\end{itemize}
so that
\begin{itemize}

  \item[] $\nextT \Skup{T} \vdash \KNI{k+1}{\Skup{\overline{B}}}{\Skup{\overline{X}}}{ \neg (( \bigwedge_{l=0}^{i-1} \nextT^l \alpha) \wedge \nextT^{i} \beta)}$, for every $i \in \Nat$

  \item[] $\nextT \Skup{T} \vdash \KNI{k+1}{\Skup{\overline{B}}}{\Skup{\overline{X}}}{ \neg (\alpha \until \beta)}$, by R$\until$

  \item[] $\nextT \Skup{T} \vdash (\beta_{k} \vee \neg \beta_{k}) \rightarrow \nextT \KNI{k}{\Skup{B}}{\Skup{X}}{ \neg (\alpha \until \beta))}$, and

  \item[] $\nextT \Skup{T} \vdash \nextT \KNI{k}{\Skup{B}}{\Skup{X}}{ \neg (\alpha \until \beta)}$.
\end{itemize}
The other cases can be proved similarly. \hfill $\blacksquare$
\end{proof}

The next lemma contains some standard statements about linear temporal, epistemic and probabilistic logics.

\begin{lemma}\label{auxlemaobjthm}
The following hold:
\begin{enumerate}
  \item if $\vdash \alpha \leftrightarrow \beta$, then $\vdash \nextT \alpha \leftrightarrow \nextT \beta$
      \label{auxlemaobjthm_1a}

  \item if $\vdash \alpha \leftrightarrow \beta$, then $\vdash \prev \alpha \leftrightarrow \prev \beta$
      \label{auxlemaobjthm_1b}

%
%
%
  \item $\vdash \know_a (\alpha \wedge \neg \alpha ) \rightarrow \know_a (\beta \wedge \neg \beta )$ \label{auxlemaobjthm_3}

  \item $\vdash \neg \prev \alpha \rightarrow \prev \neg \alpha $  \label{auxlemaobjthm_4}

%
%

  \item $\vdash \prev (\alpha \wedge \beta) \leftrightarrow (\prev \alpha \wedge \prev \beta) $ \label{auxlemaobjthm_5}

%

  \item $\vdash (\prev \alpha \vee \prev \beta) \rightarrow \prev (\alpha \vee \beta)$ \label{auxlemaobjthm_6}

%
%
%
%
%
		
  \item $\prev (\alpha \wedge \neg \alpha) \vdash  \prev \beta$ \label{auxlemaobjthm_nova_7}

  \item $\vdash (\nextT \alpha \rightarrow \nextT \beta) \leftrightarrow \nextT (\alpha \rightarrow \beta)$, \label{auxlemaobjthm_9}

  \item $\vdash (\nextT \alpha \wedge \nextT \beta) \leftrightarrow \nextT (\alpha \wedge \beta)$, \label{auxlemaobjthm_10}
		
  \item $\vdash (\nextT \alpha \vee \nextT \beta) \leftrightarrow \nextT (\alpha \vee \beta)$,\label{auxlemaobjthm_11}

%
%
		
  \item for $j \in \Nat$, $\nextT^j \beta, \nextT^0 \alpha, \ldots , \nextT^{j-1} \alpha \vdash \alpha \until \beta$ \label{auxlemaobjthm_14}

  \item $\vdash \alpha \since \beta \leftrightarrow (\prev (\alpha \wedge \neg \alpha) \wedge \beta) \vee (\neg \prev (\alpha \wedge \neg \alpha) \wedge (\beta \vee (\alpha \wedge \prev (\alpha \since \beta))))$. \label{auxlemaobjthm_15}

  \item for $j \in \Nat$, $\prev^{j} \beta, \bigwedge_{k=1}^{j} \neg \prev^{k}(\alpha \wedge \neg \alpha), \bigwedge_{l=0}^{j-1} \prev^{l}\alpha \vdash \alpha \since \beta$.
      \label{auxlemaobjthm_16}

\item $\vdash \mathtt{P}_{\geqslant 1}(\alpha \rightarrow \beta) \rightarrow (\mathtt{P}_{\geqslant s}\alpha \rightarrow \mathtt{P}_{\geqslant s} \beta)$, for every $s \in \UnitQ$. \label{auxlemaobjthm_gp19_1}

\item If $\vdash \alpha \leftrightarrow \beta$, then   $\vdash \mathtt{P}_{\geqslant s} \alpha \leftrightarrow \mathtt{P}_{\geqslant s} \beta$, for every $s \in \UnitQ$. \label{auxlemaobjthm_gp19_2}

\item $\vdash \mathtt{P}_{a, \geqslant 1}(\alpha \rightarrow \beta) \rightarrow (\mathtt{P}_{a, \geqslant s}\alpha \rightarrow \mathtt{P}_{a, \geqslant s} \beta)$, \label{auxlemaobjthm_19_1}

\item If $\vdash \alpha \leftrightarrow \beta$, then   $\vdash \mathtt{P}_{a, \geqslant s} \alpha \leftrightarrow \mathtt{P}_{a, \geqslant s} \beta$, for every $s \in \UnitQ$. \label{auxlemaobjthm_19_2}

  \item $\vdash \mathtt{P}_{\geqslant s} \alpha \rightarrow \mathtt{P}_{\geqslant r} \alpha$, for $s \geqslant r$, for every $s \in \UnitQ$. \label{auxlemaobjthm_gp17}

  \item $\vdash \mathtt{P}_{a,\geqslant s} \alpha \rightarrow \mathtt{P}_{a,\geqslant r} \alpha$, for $s \geqslant r$. \label{auxlemaobjthm_17}

  \item $\vdash \mathtt{P}_{\leqslant s} \alpha \rightarrow \mathtt{P}_{\leqslant r} \alpha$, for $r \geqslant s$. \label{auxlemaobjthm_gp18}

  \item $\vdash \mathtt{P}_{a,\leqslant s} \alpha \rightarrow \mathtt{P}_{a,\leqslant r} \alpha$, for $r \geqslant s$. \label{auxlemaobjthm_18}

\end{enumerate}
\end{lemma}


\begin{theorem}[Witnesses' theorem]\label{witnesses}
Let $\Skup{T}$ be a consistent set of formulas. Then:
\begin{enumerate}

  \item If $\Skup{T} \cup \{ \KNI{k}{\Skup{B}}{\Skup{X}}{\neg (\alpha \until \beta)} \}$ is not consistent, then there is $i_0 \in \Nat$ such that
      $$\Skup{T} \cup \{ \neg \KNI{k}{\Skup{B}}{\Skup{X}}{\neg (( \bigwedge_{l=0}^{i_0-1} \nextT^l \alpha) \wedge \nextT^{i_0} \beta)} \}$$
      is consistent. \label{witnesses_1}

  \item If $\Skup{T} \cup \{ \KNI{k}{\Skup{B}}{\Skup{X}}{\neg (\alpha \since \beta)} \}$ is not consistent, then there is $i_0 \in \Nat$ such that
      $$\Skup{T} \cup \{ \neg \KNI{k}{\Skup{B}}{\Skup{X}}{\neg ( ( \bigwedge_{l=0}^{i_0-1} \prev^l \alpha) \wedge ( \bigwedge_{l=0}^{i_0} \neg \prev^l (\alpha \wedge \neg \alpha)) \wedge \prev^{i_0} \beta)} \}$$
      is consistent. \label{witnesses_2}

  \item If $\Skup{T} \cup \{ \KNI{k}{\Skup{B}}{\Skup{X}}{\commonK \alpha} \}$ is not consistent, then there is $i_0 \in \Nat$ such that
      $$\Skup{T} \cup \{ \neg \KNI{k}{\Skup{B}}{\Skup{X}}{(\groupK)^{i_0} \alpha} \}$$
      is consistent. \label{witnesses_3}

  \item If $\Skup{T} \cup \{ \KNI{k}{\Skup{B}}{\Skup{X}}{\mathtt{P}_{\geqslant r} \alpha} \}$ is not consistent, then there is $i_0 \in \Nat$ such that
      $$\Skup{T} \cup \{ \neg\KNI{k}{\Skup{B}}{\Skup{X}}{ \mathtt{P}_{\geqslant r - \frac{1}{i_0}} \alpha} \}$$
      is consistent. \label{witnesses_4}

  \item If $\Skup{T} \cup \{ \KNI{k}{\Skup{B}}{\Skup{X}}{\mathtt{P}_{a,\geqslant r} \alpha} \}$ is not consistent, then there is $i_0 \in \Nat$ such that
      $$\Skup{T} \cup \{ \neg\KNI{k}{\Skup{B}}{\Skup{X}}{ \mathtt{P}_{a,\geqslant r - \frac{1}{i_0}} \alpha} \}$$
      is consistent. \label{witnesses_5}

\end{enumerate}
\end{theorem}
%
%
\begin{proof}
Here we consider the case (\ref{witnesses_5}):
Let $\Skup{T}$ be a consistent set of formulas such that:
\begin{itemize}

  \item $\Skup{T} \cup \{ \KNI{k}{\Skup{B}}{\Skup{X}}{\mathtt{P}_{a,\geqslant r} \alpha} \}$ is not consistent, and

  \item for every $i_0 \in \Nat$, $\Skup{T} \cup \{ \neg\KNI{k}{\Skup{B}}{\Skup{X}}{ \mathtt{P}_{a,\geqslant r - \frac{1}{i_0}} \alpha} \}$ is not consistent.
\end{itemize}
Then:
\begin{itemize}
  \item[] $\Skup{T} \cup \{ \neg \KNI{k}{\Skup{B}}{\Skup{X}}{ \mathtt{P}_{a,\geqslant r - \frac{1}{i_0}} \alpha} \} \vdash (\gamma \wedge \neg \gamma)$, for every $i_0 \in \Nat$

  \item[] $\Skup{T} \vdash \neg \KNI{k}{\Skup{B}}{\Skup{X}}{ \mathtt{P}_{a,\geqslant r - \frac{1}{i_0}} \alpha}  \rightarrow (\gamma \wedge \neg \gamma)$, by Deduction theorem, for every $i_0 \in \Nat$

  \item[] $\Skup{T} \vdash \KNI{k}{\Skup{B}}{\Skup{X}}{ \mathtt{P}_{a,\geqslant r - \frac{1}{i_0}} \alpha}$, for every $i_0 \in \Nat$

  \item[] $\Skup{T} \vdash \KNI{k}{\Skup{B}}{\Skup{X}}{\mathtt{P}_{a,\geqslant r} \alpha}$, by RA,
\end{itemize}
which contradicts the assumption about consistency of $\Skup{T}$. Hence, there is $i_0 \in \Nat$ such that
$$\Skup{T} \cup \{ \neg\KNI{k}{\Skup{B}}{\Skup{X}}{ \mathtt{P}_{a,\geqslant r - \frac{1}{i_0}} \alpha} \}$$
is consistent. \hfill $\blacksquare$
\end{proof}


The next lemma will be used in the following proofs.

\begin{lemma}\label{auxlemamcs}
Let $\Skup{T}$ be a maximal consistent consistent set of formulas. Then:
\begin{enumerate}

\item For every formula $\alpha$, exactly one of $\alpha$ and $\neg \alpha$ is in $\Skup{T}$. \label{auxlemamcs_1}

\item $\Skup{T}$ is deductively closed. \label{auxlemamcs_2}

\item $\alpha \land \beta \in \Skup{T}$ iff $\alpha \in \Skup{T}$ and $\beta \in \Skup{T}$. \label{auxlemamcs_3}

\item \label{imp} If $\{\alpha, \alpha \to \beta\} \subseteq \Skup{T}$, then $\beta \in \Skup{T}$. \label{auxlemamcs_4}

\item If $\{ \KNI{k}{\Skup{B}}{\Skup{X}}{ \neg (( \bigwedge_{l=0}^{i-1} \nextT^l \alpha) \wedge \nextT^{i} \beta)}\, : \,  i \in \Nat\} \subset \Skup{T}$, then $\KNI{k}{\Skup{B}}{\Skup{X}}{ \neg (\alpha \until \beta)} \in \Skup{T}$. \label{auxlemamcs_5}

\item If $\{ \KNI{k}{\Skup{B}}{\Skup{X}}{ \neg ( ( \bigwedge_{l=0}^{i-1} \prev^l \alpha) \wedge ( \bigwedge_{l=0}^{i} \neg \prev^l (\alpha \wedge \neg \alpha)) \wedge \prev^{i} \beta)}\, : \,  i \in \Nat\} \subset \Skup{T}$, then $\KNI{k}{\Skup{B}}{\Skup{X}}{ \neg (\alpha \since \beta)} \in \Skup{T}$. \label{auxlemamcs_6}

\item If $\{ \KNI{k}{\Skup{B}}{\Skup{X}}{(\groupK)^i \alpha}\, : \,  i \in \Nat\} \subset \Skup{T}$, then $\KNI{k}{\Skup{B}}{\Skup{X}}{\commonK \alpha} \in \Skup{T}$. \label{auxlemamcs_7}

\item If $r = \sup \,\{s  \in \UnitQ \, :  \, \mathtt{P}_{\geqslant s} \alpha \in \Skup{T} \}$, and $r \in \UnitQ$, then $\mathtt{P}_{\geqslant r} \alpha \in \Skup{T}$. \label{auxlemamcs_gp8}

\item If $r = \sup \,\{s  \in \UnitQ \, :  \, \mathtt{P}_{\geqslant s} \alpha \in \Skup{T}\}$, then for every $s \in \UnitQ$ such that $r > s$, $\mathtt{P}_{\geqslant s} \alpha \in \Skup{T}$. \label{auxlemamcs_gp9}

\item If $r = \sup \,\{s  \in \UnitQ \, :  \, \mathtt{P}_{a,\geqslant s} \alpha \in \Skup{T} \}$, and $r \in \UnitQ$, then $\mathtt{P}_{a,\geqslant r} \alpha \in \Skup{T}$. \label{auxlemamcs_8}

\item If $r = \sup \,\{s  \in \UnitQ \, :  \, \mathtt{P}_{a,\geqslant s} \alpha \in \Skup{T}\}$, then for every $s \in \UnitQ$ such that $r > s$, $\mathtt{P}_{a,\geqslant s} \alpha \in \Skup{T}$. \label{auxlemamcs_9}

\item There is a positive $i \in \Nat$ such that for every $\alpha \in \Fle$, $\prev^i (\alpha \wedge \neg \alpha) \in \Skup{T}$, and for every $j \in \Nat$, $j<i$, $\prev^j (\alpha \wedge \neg \alpha) \not\in \Skup{T}$. \label{auxlemamcs_10}

\item $\Skup{T}^{-\nextT}$ is maximal consistent. \label{auxlemamcs_11}

\item If for every $\alpha \in \Fle$, $\prev (\alpha \wedge \neg \alpha) \not\in \Skup{T}$, then $\Skup{T}^{-\prev}$ is maximal consistent. \label{auxlemamcs_12}

\item If for every $\alpha \in \Fle$, $\prev (\alpha \wedge \neg \alpha) \not\in \Skup{T}$, then $\Skup{T} = (\Skup{T}^{-\prev})^{-\nextT}$. \label{auxlemamcs_13}

\item If $\know_a \alpha \not\in \Skup{T}$, then $\Skup{T}^{-\know_a} \cup \{ \neg \alpha \}$ is consistent. \label{auxlemamcs_14} \hfill $\blacksquare$
\end{enumerate}
\end{lemma}


\begin{theorem}[Lindenbaum's theorem]\label{max_ext}
Every $\Ax$-consistent set of formulas $\Skup{T}$ can be extended to a maximal $\Ax$-consistent set $\Skup{T}^*$.
\end{theorem}
\begin{proof}
The proof will be based on the following ideas:
\begin{itemize}
  \item a procedure for extending $\Skup{T}$ will be described so that in each step a consistent superset of $\Skup{T}$ is obtained (by adding new formulas to supersets of $\Skup{T}$),

  \item the procedure guarantees that, if a formula which is the negation of a conclusion of an infinitary rule is added to a superset of $\Skup{T}$, then a witness (the negation of a premise of the rule) is also added to the superset,

  \item it will be shown that the union of all those extensions:
  \begin{itemize}
    \item contains exactly one of $\alpha$, $\neg \alpha$ for every $\alpha \in \Fle$, and
    \item is deductively closed sets,
  \end{itemize}

  \item so that the union is a maximal consistent extension of $\Skup{T}$.
\end{itemize}

Let $\{ \alpha_i \, | \,  i \in \Nat\}$ be a list of all $\Fle$-formulas. We define a sequence of theories $\{ \Skup{T}_i \, | \,  i \in \Nat\}$ and a theory $\Skup{T}^\ast$ as follows:
\begin{enumerate}
  \item \label{lin1} $\Skup{T}_0 = \Skup{T}$.

  \item \label{lin2} For every $i \in \Nat$:
  \begin{enumerate}
    \item \label{lin2_1} If $\Skup{T}_i \cup \{ \alpha_i \}$ is consistent, then $\Skup{T}_{i+1} = \Skup{T}_i \cup \{ \alpha_i \}$.

    \item \label{lin2_2} If $\Skup{T}_i \cup \{ \alpha_i \}$ is inconsistent, then:
    \begin{enumerate}
      \item \label{lin2_2_1} If $\alpha_i = \KNI{k}{\Skup{B}}{\Skup{X}}{ \neg (\alpha \until \beta)}$, then
      $$\Skup{T}_{i+1} = \Skup{T}_i \cup \{ \neg \alpha_i, \neg \KNI{k}{\Skup{B}}{\Skup{X}}{ \neg (( \bigwedge_{l=0}^{j-1} \nextT^l \alpha) \wedge \nextT^{j} \beta)}\}$$
      for some $j \in \Nat$ so that $\Skup{T}_{i+1}$ is consistent.

      \item \label{lin2_2_2} If $\alpha_i = \KNI{k}{\Skup{B}}{\Skup{X}}{ \neg (\alpha \since \beta)}$, then
      $$\Skup{T}_{i+1} = \Skup{T}_i \cup \{ \neg \alpha_i, \neg  \KNI{k}{\Skup{B}}{\Skup{X}}{ \neg ( ( \bigwedge_{l=0}^{j-1} \prev^l \alpha) \wedge ( \bigwedge_{l=0}^{j} \neg \prev^l (\alpha \wedge \neg \alpha)) \wedge \prev^{i} \beta)} \}$$
      for some $j \in \Nat$ so that $\Skup{T}_{i+1}$ is consistent.

      \item \label{lin2_2_3} If $\alpha_i =  \KNI{k}{\Skup{B}}{\Skup{X}}{\commonK \alpha}$, then
      $$\Skup{T}_{i+1} = \Skup{T}_i \cup \{ \neg \alpha_i, \neg \KNI{k}{\Skup{B}}{\Skup{X}}{(\groupK)^j \alpha} \}$$
      for some $j \in \Nat$ so that $\Skup{T}_{i+1}$ is consistent.

      \item \label{lin2_2_gp4} If $\alpha_i = \KNI{k}{\Skup{B}}{\Skup{X}}{\mathtt{P}_{\geqslant r} \alpha}$, then
      $$\Skup{T}_{i+1} = \Skup{T}_i \cup \{ \neg \alpha_i, \neg \KNI{k}{\Skup{B}}{\Skup{X}}{ \mathtt{P}_{\geqslant r - \frac{1}{j}} \alpha} \}$$
      for some $j \in \Nat$ so that $\Skup{T}_{i+1}$ is consistent.

      \item \label{lin2_2_4} If $\alpha_i = \KNI{k}{\Skup{B}}{\Skup{X}}{\mathtt{P}_{a,\geqslant r} \alpha}$, then
      $$\Skup{T}_{i+1} = \Skup{T}_i \cup \{ \neg \alpha_i, \neg \KNI{k}{\Skup{B}}{\Skup{X}}{ \mathtt{P}_{a,\geqslant r - \frac{1}{j}} \alpha} \}$$
      for some $j \in \Nat$ so that $\Skup{T}_{i+1}$ is consistent.

      \item \label{lin2_2_5} Otherwise, $\Skup{T}_{i+1} = \Skup{T}_i \cup \{ \neg \alpha_i \}$.
    \end{enumerate}
  \end{enumerate}

  \item \label{lin3} $\Skup{T}^\ast = \cup_{i \in \Nat} \Skup{T}_i$.
\end{enumerate}
First, we prove that all theories $\Skup{T}_{i}$ are consistent. Note that this trivially holds for theories obtained by the steps \ref{lin1}, \ref{lin2_1}, and \ref{lin2_2_5} of the above construction. Theorem \ref{witnesses} guarantees that the same holds for the steps \ref{lin2_2_1} -- \ref{lin2_2_4}.

Second, we show that $\Skup{T}^\ast$ is a maximal consistent set of formulas. We start by noticing that the steps \ref{lin2_1} and \ref{lin2_2} of above construction guarantee that for every $\alpha \in \Fle$, at least one of $\alpha$ and $\neg \alpha$ belongs to $\Skup{T}^\ast$. On the other hand, it is not possible that both $\alpha$ and $\neg \alpha$ are in $\Skup{T}^\ast$, since there would exist $i$ and $\Skup{T}_{i+1}$ such that:
\begin{itemize}
  \item $\alpha$ and $\neg \alpha$ are $\alpha_j$, and $\alpha_k$, respectively, from the above enumeration of all $\Fle$-formulas,
  \item $i = \max \{ j, k \}$, and
  \item $\alpha, \neg \alpha \in \Skup{T}_{i+1}$,
\end{itemize}
which contradicts consistency of $\Skup{T}_{i+1}$.


Finally, using transfinite induction on the length of a proof we show that $\Skup{T}^\ast$ is deductively closed. Let $\Skup{T}^\ast \vdash \gamma$. Recall that if $\Skup{T}^\ast \vdash \gamma$, then $\gamma$ may be:
\begin{itemize}
  \item an instance of some axiom schemata, or
  \item a formula from the set $\Skup{T}^\ast$, or
  \item obtained from the previous formulas by an inference rule (where premises of the inference rules R$\nextT$N, R$\prev$N, R$\know_a$N, RGPN and RPN must be theorems).
\end{itemize}
If $\gamma$ is an instance of an axiom schemata, then $\neg \gamma$ cannot belong to any $\Skup{T}_{k}$ since
\begin{itemize}
  \item[] $\Skup{T}_{k} \vdash \gamma$, and
  \item[] $\Skup{T}_{k} \vdash \neg \gamma$
\end{itemize}
contradict consistency of $\Skup{T}_{k}$. If $\gamma \in \Skup{T}^\ast$, the statement trivially holds.

Let  $\gamma$ be obtained from $\Skup{T}^\ast$ by an application of one of the finitary rules. By the induction hypothesis, there is $l \in \Nat$ such that all premisses of the rule belong to $\Skup{T}_{l}$. Let $\gamma = \alpha_j$ and $\neg \gamma = \alpha_k$ in the above enumeration of all $\Fle$-formulas, and $i = \max \{ j, k, l \}$. If $\neg \gamma \in \Skup{T}_{i+1}$, then:
\begin{itemize}
  \item[] $\Skup{T}_{i+1} \vdash \gamma$, and
  \item[] $\Skup{T}_{i+1} \vdash \neg \gamma$,
\end{itemize}
which contradicts consistency of $\Skup{T}_{i+1}$.


Finally, let $\gamma$ be a consequence of one of the infinitary rules R$\until$, R$\since$, RC, RGA and RA. Let us consider Rule R$\until$.


Let $\gamma = \KNI{k}{\Skup{B}}{\Skup{X}}{ \neg (\alpha \until \beta)}$ be obtained by an application of Rule R$\until$:
\begin{itemize}
  \item[] $\Skup{T}^\ast \vdash \KNI{k}{\Skup{B}}{\Skup{X}}{ \neg (( \bigwedge_{l=0}^{m-1} \nextT^l \alpha) \wedge \nextT^{m} \beta)}$, for every $m \in \Nat$ and
  \item[] $\Skup{T}^\ast \vdash \KNI{k}{\Skup{B}}{\Skup{X}}{ \neg (\alpha \until \beta)}$, by R$\until$.
\end{itemize}
By the induction hypothesis, every $\KNI{k}{\Skup{B}}{\Skup{X}}{ \neg (( \bigwedge_{l=0}^{m-1} \nextT^l \alpha) \wedge \nextT^{m} \beta)} \in \Skup{T}^\ast$.
If $\KNI{k}{\Skup{B}}{\Skup{X}}{ \neg (\alpha \until \beta)} \not\in \Skup{T}^\ast$, let $\KNI{k}{\Skup{B}}{\Skup{X}}{ \neg (\alpha \until \beta)} = \alpha_i$ in the above enumeration of all $\Fle$-formulas. By the step \ref{lin2_2_1} of the above construction, there is $j \in \Nat$ so that $\neg  \KNI{k}{\Skup{B}}{\Skup{X}}{ \neg (( \bigwedge_{l=0}^{j-1} \nextT^l \alpha) \wedge \nextT^{j} \beta)} \in \Skup{T}_{i+1}$.
Let $\KNI{k}{\Skup{B}}{\Skup{X}}{ \neg (( \bigwedge_{l=0}^{j-1} \nextT^l \alpha) \wedge \nextT^{j} \beta)} = \alpha_k$ in the above enumeration of all $\Fle$-formulas.
By the induction hypothesis $\KNI{k}{\Skup{B}}{\Skup{X}}{ \neg (( \bigwedge_{l=0}^{j-1} \nextT^l \alpha) \wedge \nextT^{j} \beta)} \in \Skup{T}^\ast$, so $\KNI{k}{\Skup{B}}{\Skup{X}}{ \neg (( \bigwedge_{l=0}^{j-1} \nextT^l \alpha) \wedge \nextT^{j} \beta)} \in \Skup{T}_{k+1}$. Then we have:
\begin{itemize}
  \item $\neg  \KNI{k}{\Skup{B}}{\Skup{X}}{ \neg (( \bigwedge_{l=0}^{j-1} \nextT^l \alpha) \wedge \nextT^{j} \beta)} \in \Skup{T}_{i+1}$, and
  \item $\KNI{k}{\Skup{B}}{\Skup{X}}{ \neg (( \bigwedge_{l=0}^{j-1} \nextT^l \alpha) \wedge \nextT^{j} \beta)} \in \Skup{T}_{k+1}$,
\end{itemize}
which contradicts consistency of $\Skup{T}_{l}$, where $l = \max \{ i+1, k+1 \}$, since $\Skup{T}_{i+1}, \Skup{T}_{k+1} \subset \Skup{T}_{l}$.

In this way we show that $\Skup{T}^\ast$ is deductively closed and does not contain all formulas, it is consistent. As it is noted above, for each $\alpha \in \Fle$ exactly one of $\alpha$ and $\neg \alpha$ belongs to $\Skup{T}^\ast$. Thus, $\Skup{T}^\ast$ is a maximal consistent set. \hfill $\blacksquare$
\end{proof}


\section{Strong completeness}

First note that the following can be showed using straightforward, but tedious proof:
\begin{theorem}[Soundness for $\Ax$]\label{soundness}
	$ \vdash \beta$ implies $ \models \beta$. \hfill $\blacksquare$
\end{theorem}

Next, we construct the canonical model in which the considered consistent set $\Skup{T}$ is satisfied. Note that Theorem \ref{max_ext} gives us a maximal consistent set $\Skup{T}^\ast$ which extends $\Skup{T}$. Then:
\begin{itemize}
  \item By Lemma \ref{auxlemamcs}.\ref{auxlemamcs_10} there is a unique positive $k \in \Nat$ such that for every $\alpha \in \Fle$, $\prev^k (\alpha \wedge \neg \alpha) \in \Skup{T}^\ast$, and for every $j \in \Nat$, $j < k$, for every $\alpha \in \Fle$, $\prev^j (\alpha \wedge \neg \alpha) \not\in \Skup{T}^\ast$.

  \item By Lemma \ref{auxlemamcs}.\ref{auxlemamcs_12} the set $\Skup{T}^{\ast^{-\prev^{k-1}}}$ is a maximal consistent set.

  \item The set $\Skup{T}^{\ast^{-\prev^{k-1}}}_{0} = \{ \prev (\alpha \wedge \neg \alpha) \, : \, \alpha \in \Fle \}$ is consistent, since it is a subset of  $\Skup{T}^{\ast^{-\prev^{k-1}}}$.
\end{itemize}
The idea of our construction is that possible worlds of the canonical model correspond to maximal consistent sets of formulas, while we define sequences of maximal consistent sets and the corresponding runs as follows:
\begin{itemize}
  \item Let the set $\Skup{Mcs}$ contain all maximal consistent extensions of the set $\Skup{T}^{\ast^{-\prev^{k-1}}}_{0}$. Elements of $\Skup{Mcs}$ correspond to starting possible worlds of runs in the canonical model. Note that Theorem \ref{max_ext} guarantees that $\Skup{Mcs}$ is nonempty.

  \item For an arbitrary $\Skup{S} \in \Skup{Mcs}$, we define the sequence of maximal consistent sets and the corresponding run $r^{\Skup{S}}$ in the following way:
\begin{itemize}
  \item $\Skup{S}_0 = \Skup{S}$.

  \item For $i \in \Nat$, $\Skup{S}_{i+1} = \Skup{S}_{i}^{-\nextT}$.

  \item By Lemma \ref{auxlemamcs}.\ref{auxlemamcs_11} every $\Skup{S}_{i+1}$ is maximal consistent.

  \item For $i \in \Nat$ and every propositional letter $p\in \SetOfPropLet$, $p \in r^{\Skup{S}}(i)$ iff $p \in \Skup{S}_{i}$.

\end{itemize}

  \item Note that Lemma \ref{auxlemamcs}.\ref{auxlemamcs_13} guaranties that there is a sequence of maximal consistent sets $\Skup{S}_0$, $\Skup{S}_1$, \ldots, such that $\Skup{S}_{k-1} = \Skup{T}^\ast$.

  \item Note that Lemma \ref{auxlemaobjthm}(\ref{auxlemaobjthm_nova_7}) guaranties that for every $\alpha \in \Fle$, and every $\Skup{S} \in \Skup{Mcs}$, $\prev \alpha \in \Skup{S}$.
\end{itemize}
Let $\Skup{R}^{\ast}$ contain all such runs and $\Skup{W}^{\ast}$ be the set of all possible worlds. Then we define:
\begin{itemize}

\item $\mathcal{A}^{\ast}: \Skup{W}^{\ast} \rightarrow \Powerset{\SetOfAgents}$ such that for every agent $a \in \SetOfAgents$:
    \begin{itemize}
    \item $a \in \mathcal{A}^{\ast}((r^{\Skup{S}},n))$ iff $A_a \in \Skup{S}_{n}$.
    \end{itemize}

\item $\mathcal{K}^{\ast}$, the set of accessibility relations on possible worlds such that for every agent $a \in \SetOfAgents$:
\begin{itemize}
  \item If $a \not\in \mathcal{A}^{\ast}((r^{\Skup{S}},n))$, then $(r^{\Skup{S}},n) \mathcal{K}^{\ast}_a (r',n')$ is false for all $(r',n')$, otherwise

  \item $(r^{\Skup{S}},n) \mathcal{K}^{\ast}_a (r^{\Skup{S'}},n')$ iff $\Skup{S}_n^{-\know_a} = \{ \alpha: \know_a \alpha \in \Skup{S}_n \} \subset \Skup{S'}_{n'}$.
\end{itemize}

\item $\mathcal{P}^{\ast}$ is a function on the set $\Skup{W}^{\ast}$ such that $\mathcal{P}^{\ast}((r^{\Skup{S}},n)) = \langle H^{\ast,(r^{\Skup{S}},n)}, \mu^{\ast,(r^{\Skup{S}},n)}, \{\mathcal{P}^{\ast}_a: a \in \SetOfAgents\}\rangle$ and:
\begin{itemize}

\item for every $\alpha \in \Fle$, $[\![\alpha]\!]^{(r^{\Skup{S}},n)} = \{ r^{\Skup{S'}} \in \Skup{R}^{\ast} \, : \, \alpha\in \Skup{S'}_{0} \}$,

\item $H^{\ast,(r^{\Skup{S}},n)}$ is a family of sets $\{ [\![\alpha]\!]^{(r^{\Skup{S}},n)} \, : \, \alpha \in \Fle\}$,

\item $\mu^{\ast,(r^{\Skup{S}},n)}: H^{\ast,(r^{\Skup{S}},n)} \rightarrow [0,1]$ such that $\mu^{\ast,(r^{\Skup{S}},n)}([\![\alpha]\!]^{(r^{\Skup{S}},n)}) = \sup \{ t \in \UnitQ \, : \, P_{\geqslant t} \alpha \in \Skup{S}_{n} \}$, and

\item $\mathcal{P}^{\ast}$ is a set of functions $\mathcal{P}^{\ast}_a$ for every agent $a \in \SetOfAgents$, where for every possible world $(r^{\Skup{S}},n) \in \Skup{W}^{\ast}$, $\mathcal{P}^{\ast}_a((r^{\Skup{S}},n)) = \langle \Skup{W}^{\ast}_a(r^{\Skup{S}},n), H^{\ast}_a(r^{\Skup{S}},n), \mu^{\ast}_a(r^{\Skup{S}},n) \rangle$ and:
\begin{itemize}
  \item $\Skup{W}^{\ast, (r^{\Skup{S}},n)}_a = \Skup{W}^{\ast}$,

  \item for every $\alpha \in \Fle$, $[\![\alpha]\!]_a^{(r^{\Skup{S}},n)} = \{ (r^{\Skup{S'}},n') \in \Skup{W}^{\ast, (r^{\Skup{S}},n)}_a \, : \, \alpha\in \Skup{S'}_{n'} \}$,

  \item $H^{\ast, (r^{\Skup{S}},n)}_a$ is a family of sets $\{ [\![\alpha]\!]_a^{(r^{\Skup{S}},n)} \, : \, \alpha \in \Fle\}$, and

  \item $\mu^{\ast, (r^{\Skup{S}},n)}_a: H_a(r^{\Skup{S}},n) \rightarrow [0,1]$ such that $\mu^{\ast, (r^{\Skup{S}},n)}_a([\![\alpha]\!]_a^{(r^{\Skup{S}},n)}) = \sup \{ t \in \UnitQ \, : \, P_{a, \geqslant t} \alpha \in \Skup{S}_n \}$.
\end{itemize}
\end{itemize}
\end{itemize}

\begin{definition} \label{canonicak_model_def}
Let $\Skup{T}$ be a consistent set and $\Skup{T}^\ast$ its maximal consistent extension. Let $\Skup{R}^{\ast}$, $\mathcal{A}^{\ast}$, $\mathcal{K}^{\ast}$ and $\mathcal{P}^{\ast}$ be defined as above. The canonical model which corresponds to the set $\Skup{T}$ is $\Model{M}^{\ast} = \langle \Skup{R}^{\ast}, \mathcal{A}^{\ast}, \mathcal{K}^{\ast}, \mathcal{P}^{\ast} \rangle$.
\hfill $\blacksquare$
\end{definition}

Note that Definition \ref{canonicak_model_def} of the canonical model $\Model{M}^{\ast}$ relies on sets of the forms $[\![\alpha]\!]^{(r^{\Skup{S}},n)}$ and $[\![\alpha]\!]_a^{(r^{\Skup{S}},n)}$, while in Definition \ref{def_modeli} sets of the forms $[\alpha]^{(r,n)}$ and $[\alpha]_a^{(r,n)}$ are used. Recall that $[\alpha]^{(r,n)}$ and $[\alpha]_a^{(r,n)}$ are defined using the satisfiability relation, which is not the case for $[\![\alpha]\!]^{(r^{\Skup{S}},n)}$ and $[\![\alpha]\!]_a^{(r^{\Skup{S}},n)}$. So, having the notion of the canonical model $\Model{M}^{\ast}$ we have first to prove that $\Model{M}^{\ast}$ is indeed a $\ClassModels$-model, where the hardest part is to show that the sets $[\![\alpha]\!]^{(r^{\Skup{S}},n)}$ and $[\alpha]^{(r^{\Skup{S}},n)}$ (i.e., $[\![\alpha]\!]_a^{(r^{\Skup{S}},n)}$ and $[\alpha]_a^{(r^{\Skup{S}},n)}$) coincide. We have to show that the families $H^{\ast,(r^{\Skup{S}},n)}$ and $H^{\ast, (r^{\Skup{S}},n)}_a$ are algebras, while the mappings $\mu^{\ast,(r^{\Skup{S}},n)}$ and $\mu^{\ast, (r^{\Skup{S}},n)}_a$ are finitely-additive probability measures. It can be proved using the ideas from \cite{ORM16}:

\begin{lemma} \label{canonical_model_meausure}
Let the canonical model $\Model{M}^{\ast}$ be defined as above. Then for every possible world $(r^{\Skup{S}},n)$ and for every agent $a \in \SetOfAgents$:
\begin{enumerate}
  \item $H^{\ast,(r^{\Skup{S}},n)}$ is an algebra of sets, \label{canonical_model_meausure_1}
  \item $\mu^{\ast,(r^{\Skup{S}},n)}$ is a finitely-additive probability measure, \label{canonical_model_meausure_2}
  \item $H^{\ast, (r^{\Skup{S}},n)}_a$ is an algebra of sets, and \label{canonical_model_meausure_3}
  \item $\mu^{\ast,(r^{\Skup{S}},n)}_a$ is a finitely-additive probability measure. \label{canonical_model_meausure_4}  \hfill $\blacksquare$
\end{enumerate}
\end{lemma}

Having Lemma \ref{canonical_model_meausure} we ensure that:

\begin{theorem} \label{canonical_model_th}
The canonical model $\Model{M}^{\ast}$ is a $\ClassModels$-model. \hfill $\blacksquare$
\end{theorem}

Finally, strong completeness follows from the previous statements.

\begin{theorem}[Strong completeness for $\Ax$]\label{completness}
	A set $\Skup{T}$ of formulas is $\Ax$-consistent iff it is satisfiable.
\end{theorem}
\begin{proof}
The ($\Leftarrow$)-direction follows from the soundness of the above axiomatic system.
To prove the ($\Rightarrow$)-direction assume that $\Skup{T}$ is consistent. Following the above proven statements, we see that Theorem \ref{max_ext} guarantees that $\Skup{T}$ can be extended to a maximal consistent $\Skup{T}^*$, and that Theorem \ref{canonical_model_th} shows that the canonical model $\Model{M}^{\ast}$ can be defined using $\Skup{T}^*$ so that $\Skup{T}$ is satisfiable in a possible world from $\Model{M}^{\ast}$.
\hfill $\blacksquare$
\end{proof}


\section*{Acknowledgment}
This work was supported by the Science Fund of the Republic of Serbia, Grant AI4TrustBC: Advanced Artificial Intelligence Techniques for Analysis and Design of System Components Based on Trustworthy BlockChain Technology and by  the Serbian Ministry of Education, Science and Technological Development (Agreement No. 451-03-9/2021-14/200029).


\bibliographystyle{plain}
\bibliography{D1.1-AI4TrustBC}                


\end{document}